\documentclass{article}

\usepackage{latexsym,amsmath,amsfonts,amsthm,amssymb,graphicx}

\numberwithin{equation}{section}

\newtheorem{Th}{Theorem}[section]

\newtheorem{Lemma}[Th]{Lemma}

\theoremstyle{definition}
\newtheorem{Example}[Th]{Example}
\newtheorem{df}[Th]{Definition}

\theoremstyle{remark}
\newtheorem{Remark}[Th]{Remark}

\begin{document}

\title{Nonessential Functionals in Multiobjective
Optimal Control Problems\footnote{Research Report CM06/I-33.
Presented at the 5th Junior European Meeting
on Control \& Information Technology (JEM'06),
September 20--22, 2006, Tallinn, Estonia.}}

\author{Agnieszka B. Malinowska${}^{\dag}$\\
        \texttt{abmalina@pb.bialystok.pl}
        \and
        Delfim F. M. Torres${}^{\ddag}$\\
        \texttt{delfim@mat.ua.pt}}

\date{${}^{\dag}$Institute of Mathematics and Physics\\
      Technical University of Bia\l ystok \\
      15-351 Bia\l ystok, Poland\\[0.3cm]
      ${}^{\ddag}$Department of Mathematics\\
      University of Aveiro\\
      3810-193 Aveiro, Portugal}

\maketitle


\begin{abstract}
We address the problem of obtaining well-defined criteria for
multiobjective optimal control systems. Necessary and sufficient
conditions for an optimal control functional to be nonessential
are proved. The results provide effective tools for determining
nonessential objectives in vector-valued optimal control problems.
\end{abstract}

\smallskip

\noindent \textbf{Mathematics Subject Classification 2000:} 49K15, 49N90, 90C29.

\smallskip


\smallskip

\noindent \textbf{Keywords:} multiobjective dynamic optimization,
multiobjective control, Pareto-optimal control,
essential/nonessential functionals.


\section{Introduction}

Multiobjective optimal control attracts more and more attention
and is source of strong current research (see \textrm{e.g.}
\cite{AubinFrankowska,Hamel,Torres} and references therein). We
consider multiobjective problems of optimal control governed by
ordinary differential dynamical systems. This comprises an
important class of problems which naturally appear on practical
applications to Economic \cite[Chap.~8]{Mordukhovich} and
Engineering modelling \cite{Salukvadze}. Our main goal is to
extend the results found in the literature on nonessential
functions of mathematical static optimization programming
\cite{GalHanne,Malinowska2006} to functionals of optimal control
theory.

It is well known that the concept of Pareto optimality or
efficiency play a crucial role on optimal control
\cite{Leitmann,Salukvadze}. The question of obtaining well-defined
criteria for multiple criteria decision making problems seems,
however, being considered in the literature only for static
multiobjective optimization problems (\textrm{cf.}
\cite{GalHanne,Malinowska2006} and references therein). In this
work we investigate the problem of obtaining well-defined criteria
for multicriteria optimal control dynamical systems.

One of the approaches dealing with the problem of obtaining
well-defined criteria for multiple criteria static decision making
problems is the concept of nonessential objective functions. A
certain objective function is called nonessential if it does not
influence the set of efficient solutions of the vector-valued
optimization problem, that is, the set of efficient solutions is
the same both with or without that objective function. Information
about nonessential objectives helps a decision maker to know and
to understand better the problem and this might be a good starting
point for further investigation or revision of the mathematical
model. Dropping nonessential functions leads to a problem with a
smaller number of objectives, which can be solved more easily. For
this reason, the issue of nonessential objectives is a substantial
feature for multiple criteria decision making
\cite{GalHanne,GalLeberling,Malinowska2002}. To the best of the
authors knowledge, no study has been done in this field for
optimal control problems. We are interested in generalizing the
previous results on nonessential objectives found in the
literature to cover optimal control problems with a vector-valued
functional to minimize. More precisely, we generalize the concept
of nonessential objective to multicriteria functionals of optimal
control systems and we give the first steps on the corresponding
theory. Main results provide methods for identifying nonessential
objectives in nonlinear and optimal control vector-valued
optimization problems.


\section{Optimal control with a vector-valued cost}

We consider a dynamical control system described by $n$
\emph{state variables} $x = \left(x_1,\ldots,x_n\right) \in
\mathbb{R}^n$ and $r$ \emph{control variables} $u =
\left(u_1,\ldots,u_r\right) \in \mathbb{R}^r$, $r \le n$. Both
state and control variables vary with respect to the scalar
variable $t \in \mathbb{R}$. Given a control vector function
$u : [a,b] \rightarrow \mathbb{R}^r$,
the state evolution over $[a,b]$, namely $x : [a,b] \rightarrow \mathbb{R}^n$,
must satisfy the control system
\begin{equation}
\label{2}
\dot{x}(t)=h(t,x(t),u(t)) \, ,
\end{equation}
the boundary conditions
\begin{equation}
\label{1}
x(a)=\alpha\, , \quad x(b)=\beta \, ,
\end{equation}
and $m$ inequality constraints
\begin{equation}
\label{3}
g_{i}(t,x(t),u(t)) \leq 0\, , \quad i=1,\ldots m \, .
\end{equation}
We would like to find a piecewise-continuous control
function $u(\cdot)$ and the corresponding state trajectory $x(\cdot)$, satisfying
\eqref{2}, \eqref{1} and \eqref{3},
which minimizes a finite number $N$ of cost functionals,
called the optimal control multiobjective criteria or an
optimal control multiobjective performance-index:
\begin{equation*}
\min \int^{b}_{a}f(t,x(t),u(t))dt
= \min \left(\int^{b}_{a}f_{1}(t,x(t),u(t))dt,\cdots,
\int^{b}_{a}f_{N}(t,x(t),u(t))dt\right) \, .
\end{equation*}
All functions $f(t,x,u)$, $g(t,x,u)$
and $h(t,x,u)$ are assumed to be continuously differentiable
with respect to $t$ and $x$ variables. To simplify notation, we write
\begin{equation*}
    I^{N}[x,u]=\int^{b}_{a}f(t,x(t),u(t))dt
\end{equation*}
and
\begin{equation*}
    I_{i}[x,u]=\int^{b}_{a}f_{i}(t,x(t),u(t))dt\, , \quad i=1,\ldots,N \, .
\end{equation*}
In general does not exist a pair of functions $\left(x,u\right)$
that renders the minimum value to each cost functional $I_i$, $i =
1,\ldots,N$, simultaneously, and one uses the concept of
Pareto-optimality. Let us denote by $S$ the set of feasible
solutions, \textrm{i.e.} the set of all admissible functions
$\left(x,u\right)$. The multiobjective control problem consists to
find all feasible solutions that are efficient in the sense of
Definition~\ref{d1}. This problem is denoted in the sequel by
$(P)$. We remark that many practical applications that appear in
engineering and economics can be written in the form of problem
$(P)$ \cite{Salukvadze}.

\begin{df}[Pareto-optimality]
\label{d1}
A pair of functions $(\tilde{x},\tilde{u})\in S$ is said to be an efficient
(Pareto-optimal) solution of the problem $(P)$ if, and only if,
there exists no $(x,u) \in S$ such that $I^{N}[x,u]\leqq
I^{N}[\tilde{x},\tilde{u}]$, where
\begin{multline*}
I^{N}[x,u]\leqq I^{N}[\tilde{x},\tilde{u}] \\
\Leftrightarrow
\forall i\in\{1,\ldots,N\} : I_{i}[x,u]\leqslant
I_{i}[\tilde{x},\tilde{u}]\wedge \exists j \in\{1,\ldots,N\} :
I_{j}[x,u] < I_{j}[\tilde{x},\tilde{u}] \, .
\end{multline*}
The set of efficient solutions of $(P)$ is denoted by $S^{N}_{E}$.
\end{df}

The central result in optimal control theory is given
by the celebrated Pontryagin maximum principle \cite{MR29:3316b},
which is a necessary optimality condition. A version of the Pontryagin
maximum principle for Pareto-solutions of control systems
with multiple criteria was proved already in the sixties \cite{Chang}.
Roughly speaking, one can say that the necessary and sufficient conditions
for Pareto-optimality are obtained converting
the vector performance optimal control problem
to a family of scalar-index optimal control problems
by forming an auxiliary scalar integral functional as a function
of the vector-index and a vector of weighting parameters \cite{Liu,ReidCitron}.
For a gentle introduction to optimal control,
including necessary and sufficient conditions
and the question of existence, we refer the readers
to \cite{MackiStrauss,Pedregal} (scalar case)
and \cite{Leitmann,Salukvadze} (Pareto optimal control).
Here we just recall three basic lemmas (\textrm{cf.} \cite[Chap.~17]{Leitmann})
that relate the Pareto-solution of a multiobjective control problem
with the solutions of an appropriate scalar-valued cost problem.

\begin{Lemma}
\label{Le:0}
If the feasible pair $(\tilde{x},\tilde{u}) \in S$ is efficient for $(P)$,
then it is optimal for the scalar-valued cost
\begin{equation*}
I_i[x,u] \, , \quad i \in \{1,\ldots,N\}
\end{equation*}
subject to the constraints $(x,u) \in S$ and
\begin{equation*}
I_j[x,u] - I_j[\tilde{x},\tilde{u}] \le 0 \, , \quad
j = 1,\ldots,N \text{ and } j \ne i \, .
\end{equation*}
\end{Lemma}

Lemma~\ref{Le:0} is very useful because it implies
that the necessary conditions \cite{MR1897883,MR29:3316b}
are also necessary for Pareto-optimality in the optimal control
problem with a vector-valued cost. As with the necessary conditions,
next two lemmas reduce the sufficient conditions for Pareto-optimality
to sufficient conditions with a scalar-valued cost functional.

\begin{Lemma}
\label{Le:1}
A feasible pair $(\tilde{x},\tilde{u})\in S$ is efficient for $(P)$
if there exists a constant $\gamma \in \mathbb{R}^N$, with
$\gamma_i > 0$ for $i = 1,\ldots,N$ and $\sum_{i=1}^{N} \gamma_i = 1$, such that
\begin{equation*}
\sum_{i=1}^{N} \gamma_i I_i[x,u] \ge \sum_{i=1}^{N} \gamma_i I_i[\tilde{x},\tilde{u}]
\end{equation*}
for every $(x,u) \in S$.
\end{Lemma}

Unlike Lemma~\ref{Le:1}, not all components of $\gamma$ in the next Lemma~\ref{Le:2}
need to be nonzero. However, in Lemma~\ref{Le:2} the minimum
of $\sum_{i=1}^{N} \gamma_i I_i[x,u]$ must be achieved by a unique
$(\tilde{x},\tilde{u})\in S$.

\begin{Lemma}
\label{Le:2}
A feasible pair $(\tilde{x},\tilde{u})\in S$ is efficient for $(P)$
if there exists a constant $\gamma \in \mathbb{R}^N$, with
$\gamma_i \ge 0$ for $i = 1,\ldots,N$ and $\sum_{i=1}^{N} \gamma_i = 1$, such that
\begin{equation*}
\sum_{i=1}^{N} \gamma_i I_i[x,u] > \sum_{i=1}^{N} \gamma_i I_i[\tilde{x},\tilde{u}]
\end{equation*}
for every $(x,u) \in S$, $(x,u) \ne (\tilde{x},\tilde{u})$.
\end{Lemma}

Together with the Pontryagin maximum principle \cite{MR1897883,MR29:3316b},
Lemmas~\ref{Le:0}, \ref{Le:1} and \ref{Le:2} provide expedient tools to study
concrete multiobjective problems of optimal control (\textrm{cf.} \S\ref{sec:IE}).


\section{Nonessential functionals: main results}

We form a new  multiobjective control problem $(\tilde{P})$ from
$(P)$ by adding a new functional
$I_{N+1}[x,u]=\int^{b}_{a}f_{N+1}(t,x(t),u(t))dt$ to problem
$(P)$. Let $S^{N+1}_{E}$ denote the set of efficient solutions of the
problem $(\tilde{P})$. With this notation we introduce the
definition of nonessential functional.

\begin{df}
\label{d2}
The functional $I_{N+1}$ is said to be nonessential in $(\tilde{P})$
if, and only if, $S^{N}_{E}=S^{N+1}_{E}$. A functional which
is not nonessential will be called essential.
\end{df}

We are interested in characterizing the functionals which do not
change the set of efficient solutions (nonessential objective
functionals). Along the text we denote by $S_{i},
i=1,2,\ldots,N,N+1$, the set of solutions of the scalar optimal
control problem
\begin{equation*}
    \min I_{i}[x,u]
\end{equation*}
subject to $S$. We start with a simple example.

\begin{Example}
\label{e1} Consider a system characterized by a single state and
control variable ($n = r = 1$) that evolves according to the state
equation
\begin{equation*}
\dot{x}(t)=u(t)
\end{equation*}
with control constraint set
\begin{equation*}
U=\{u : [a,b] \rightarrow \mathbb{R} : |u(t)|\leq 1 \}.
\end{equation*}
The system is to be transferred from a given initial state
$x(0)=\xi \neq 0$ to a given terminal state $x(T)=0$ within an
unspecified bounded interval $[0,T]$.
 Functionals to be minimized are
\begin{equation*}
I_{1}=\int_{0}^{T} dt\, , \quad I_{2}=\int_{0}^{T} |u(t)|dt\, . \
\end{equation*}
Applying the Pontryagin maximum principle \cite{MR29:3316b} we
obtain:
\begin{equation*}
S_{1}=\{(x(t),u(t)): u(t)=-sgn\{\xi\}\}\, , \
 min\int_{0}^{T}dt=|\xi| \, .
\end{equation*}
and
\begin{equation*}
S_{2}=\{(x(t),u(t)): u(t)=-sgn\{\xi\}v(t)\},
\end{equation*}
where
\begin{equation*}
v(t)\in V=\{v(t): 0\leq v(t)\leq 1, t\in [0,T], v(t)\not\equiv 0
\}\, , \  min\int_{0}^{T}|-sgn\{\xi\}v(t)|dt=|\xi|
\end{equation*}
Details can be found in \cite{AthansFalalb}. It is easy to see that
$S_{1}\cap S_{2}= S_{1}$ (we can take $v(t)=1$, $t\in [0,T]$). In
this problem we have: $S_{1}=S_{E}^{1}=S_{E}^{2}\subset S_{2}$. Hence
$I_{2}$ is nonessential, but $I_{1}$ is essential (in order to see
this we need only to change indices).
\end{Example}

\begin{Lemma}
\label{t1}
One has $S^{N}_{E}\subset S^{N+1}_{E}$ if, and only if, for every
$(x,u)\in S^{N}_{E}$ the following condition holds:
\begin{equation*}
\exists (x',u')\in S:I^{N}[x',u']=I^{N}[x,u]\Rightarrow I_{N+1}[x',u']=I_{N+1}[x,u] \, .
\end{equation*}
\end{Lemma}

\begin{proof}
Let $S^{N}_{E}\subset S^{N+1}_{E}$ and assume,
on the contrary, that exists $(\tilde{x},\tilde{u})\in S^{N}_{E}$
such that
\begin{equation}\label{4}
\exists (x',u')\in S:I^{N}[x',u']=I^{N}[\tilde{x},\tilde{u}]
\end{equation}
and
\begin{equation}\label{5}
I_{N+1}[x',u']\neq I_{N+1}[\tilde{x},\tilde{u}].
\end{equation}
We conclude from (\ref{4}) that $(x',u')\in S^{N}_{E}$. Therefore
$(x',u')$ is not in $S^{N+1}_{E}$ or $(\tilde{x},\tilde{u})$ is not
in $S^{N+1}_{E}$ by (\ref{5}). This contradicts the fact that
$S^{N}_{E}\subset S^{N+1}_{E}$. Let us prove now the second implication.
If $S^{N}_{E}=\emptyset$, then $S^{N}_{E}\subset S^{N+1}_{E}$ . Let
$S^{N}_{E}\neq \emptyset$. Suppose that for every $(x,u)\in
S^{N}_{E}$ holds:
\begin{equation}\label{6}
\exists (x',u')\in S:I^{N}[x',u']=I^{N}[x,u]\Rightarrow
I_{N+1}[x',u']=I_{N+1}[x,u]
\end{equation}
and $S^{N}_{E}$ is not contained in $S^{N+1}_{E}$. In this case
there exists $(\tilde{x},\tilde{u})$ in $S^{N}_{E}$ which is not in
$S^{N+1}_{E}$. Hence
\begin{equation}
\label{7}
    \exists (\hat{x},\hat{u})\in S:I^{N+1}[\hat{x},\hat{u}] \leqq
    I^{N+1}[\tilde{x},\tilde{u}] \, .
\end{equation}
This gives $I^{N}[\hat{x},\hat{u}] =
    I^{N}[\tilde{x},\tilde{u}]$ and from (\ref{6}) we have
    $I_{N+1}[\hat{x},\hat{u}]=I_{N+1}[\tilde{x},\tilde{u}].$
Consequently $I^{N+1}[\hat{x},\hat{u}]
= I^{N+1}[\tilde{x},\tilde{u}]$, contrary to (\ref{7}).
\end{proof}

\begin{Remark}
Notice that in Example \ref{e1} the scalar optimal control problem
\begin{equation*}
    \min I_{1}[x,u]
\end{equation*}
subject to $S$ has a unique solution. Therefore, Lemma \ref{t1}
holds true for the example.
\end{Remark}

\begin{df}
\label{d3} A function $f : \mathbb{R}^N \rightarrow \mathbb{R}$ is
nondecreasing if for $y^{1}$ and $y^{2} \in \mathbb{R}^N$:  $y^{1}
\leqq y^{2}$ imply $f(y^{1}) \leq f(y^{2})$.
\end{df}

\begin{Th}
\label{t2}
 If $I_{N+1}[x,u]= f(I^{N}[x,u])$, where
 $f : \mathbb{R}^N \rightarrow \mathbb{R}$, then
 $S^{N}_{E}\subset S^{N+1}_{E}$. Furthermore,
 $S^{N}_{E}= S^{N+1}_{E}$ if function $f$ is nondecreasing on the
 set $I^{N}(S)$.
\end{Th}

\begin{proof}
Let $(x,u)\in S^{N}_{E}$. If there exists $(x',u')\in S$ such that
$I^{N}[x',u']=I^{N}[x,u]$, then $f(I^{N}[x',u'])=f(I^{N}[x,u])$ and
so $I_{N+1}[x',u']=I_{N+1}[x,u]$. Therefore $S^{N}_{E}\subset
S^{N+1}_{E}$ by Lemma \ref{t1}.

We will now show the inclusion
$S^{N+1}_{E}\subset S^{N}_{E}$. Let $(x,u)\in S$, $(x,u)$
being not an element of the set $S^{N}_{E}$. In this
 case there exists $(x',u')\in S$ such that $I^{N}[x',u']\leqq
 I^{N}[x,u]$. If $f$ is nondecreasing on $I^{N}(S)$, we know that
$I_{N+1}[x',u']=f(I^{N}[x',u'])\leq f(I^{N}[x,u])=I_{N+1}[x,u]$.
Hence $(x,u)$ is not an element of the set $S^{N+1}_{E}$.
\end{proof}

\begin{Remark}
Example \ref{e1} shows that sufficient condition in Theorem \ref{t2},
for an optimal control functional to be nonessential, is not
necessary.
\end{Remark}

\begin{Th}
\label{t3}
 Let $S_{N+1}=\{(x^{0},u^{0})\}$. If the functional
 $I_{N+1}$ is nonessential, then
 $(x^{0},u^{0})\in S^{N}_{E}$.
\end{Th}

\begin{proof}
Let $S^{N}_{E}=S^{N+1}_{E}$ . If
$(x^{0},u^{0})$ is not an element of the set $S^{N}_{E}$, then
$(x^{0},u^{0})$ is also not an element of the set $S^{N+1}_{E}$. In this
case, there exists $(x',u')\in S$ such that $I^{N+1}[x',u']\leqq
 I^{N+1}[x^{0},u^{0}]$. So $I_{N+1}[x',u']\leq I_{N+1}[x^{0},u^{0}]$.
 This is a contradiction to the assumption that
 $S_{N+1}=\{(x^{0},u^{0})\}$.
\end{proof}

 \begin{Th}
 \label{t4}
 Let the set $S$ be compact. If the functional $I_{N+1}$ is nonessential, then
 $S_{N+1}\cap S^{N}_{E}\neq \emptyset$.
\end{Th}

\begin{proof}
Consider the problem
\begin{equation}
\label{10}
    \min \int^{b}_{a}f(t,x(t),u(t))dt
\end{equation}
subject to $S_{N+1}$. Let $\tilde{S}$ denote the set of efficient
solutions of the problem (\ref{10}). By the compactness of the set
$S$, the set $\tilde{S}$ is nonempty. Let $(x^{0},u^{0})\in
\tilde{S}$. If $(x^{0},u^{0})$ is not an element of $S^{N}_{E}$,
then by assumption $(x^{0},u^{0})$ is not an element of
$S^{N+1}_{E}$. In this case, there exists $(x',u')\in S$ such that
$I^{N+1}[x',u']\leqq I^{N+1}[x^{0},u^{0}]$. Hence $(x',u')\in
S_{N+1}$. This contradicts the fact that $(x^{0},u^{0})$ is an
efficient solution of the problem (\ref{10}).
\end{proof}

\begin{Remark}
Notice that $I_{2}$ is nonessential in Example \ref{e1} and we have
$S^{1}_{E}\cap S_{2}\neq \emptyset$.
\end{Remark}

Next section provides an example of application of the obtained
results to check whether a functional is nonessential.


\section{An illustrative example}
\label{sec:IE}

We illustrate the obtained results with a
multiobjective control problem
borrowed from \cite[\S 4.3]{Salukvadze},
where $N = 3$, $n = 2$, $r = 1$, $m = 4$, $a = 0$, $b = T$, with
$T$ not fixed. We consider a mobile rocket car with mass one
running on rails on a closed region $-3 \le x_1 \le 3$
(we denote the position of the center of the car at time $t$ by $x_1(t)$),
whose movement we can control with its accelerator $u$,
where the maximum allowable acceleration is $1$ and the maximum break
power is $-1$, \textrm{i.e.}, $-1 \le u \le 1$ (negative force means
break, positive force means acceleration). The dynamics of the system
is given by Newton's second law, force equals mass times acceleration,
which in our setting reads as $u(t) = \ddot{x}_1(t)$. The problem is to move
the car from a given location to a pre-assigned
destination. If the car is at a position $x_1 = 1$ at time $t = 0$,
with no velocity, that is $\dot{x}_1(0) = 0$, we want to find a
piecewise constant function $u(t)$
that drives the car to $x_1(T) = 0$ at some instant $T > 0$.
The state of the system is given by the position $x_1(t)$
and the velocity $x_2(t) = \dot{x}_1(t)$ (where we are and
how fast we are going at each instant of time $t$).
Different cost criteria can be considered,
for example, minimizing the time $T$ (functional $I_1$ below);
maximizing the velocity at $T$ (maximizing $x_2(T)$,
which corresponds to functional $I_2$ below); and a linear
combination $I_3$ of these functionals: minimize
\begin{equation*}
I_{1} = \int_{0}^{T} 1 dt \, , \quad
I_{2} = \int_{0}^{T} - u(t) dt \, , \quad
I_{3} = I_{1} + I_{2} \, ,
 \end{equation*}
subject to the control system
\begin{equation}
\label{ex:cs}
\begin{cases}
\dot{x}_{1}=x_{2} \, ,\\
\dot{x}_{2}=u \, ;
\end{cases}
\end{equation}
to the boundary conditions
\begin{equation}
\label{ex:bc}
x_{1}(0)=1 \, , \quad x_{1}(T)=0 \, , \quad x_{2}(0)=0 \, ;
\end{equation}
and inequality constraints
\begin{equation}
\label{ex:ic}
|u| \leq 1 \, , \quad |x_{1}| \leq 3 \, .
\end{equation}
Denoting by $S_{i}$, $i = 1,2,3$, the solution set of the scalar
optimal control problem $\min I_i[x,u]$ subject to
\eqref{ex:cs}-\eqref{ex:ic}, we have:\footnote{The solutions to the
scalar optimal control problems $I_i[x,u] \rightarrow \min$
are found by application of the Pontryagin maximum principle
\cite{MR29:3316b}. Details can be found in \cite[\S 4.3]{Salukvadze}.}
$S_{1}=\{(x^{1},u^{1})\}$ with
\begin{gather*}
   u^{1}=-1 \, ,\\
   x^{1}_{1}=-\frac{t^{2}}{2}+1 \, , \quad
   x^{1}_{2}=-t \, , \\
   0 \leq t \leq T = \sqrt{2} \, ;
\end{gather*}
$S_{2}=\{(x^{2},u^{2})\}$ with
\begin{gather*}
u^{2}=
\begin{cases}
-1 & \text{ if } 0\leq t\leq 2 \, , \\
+1 & \text{ if } 2\leq t\leq T = 4+\sqrt{6} \, ,
\end{cases}\\
x^{2}_{1}=
\begin{cases}
-\frac{t^{2}}{2}+1 & \text{ if } 0 \leq t \leq 2 \, , \\
\frac{t^{2}}{2}-4t+5 & \text{ if }  2 \leq t \leq 4 +\sqrt{6} \, ,
\end{cases}
\quad
x^{2}_{2} =
\begin{cases}
-t & \text{ if } 0\leq t\leq 2 \, ,\\
t-4 & \text{ if } 2 \leq t \leq 4+\sqrt{6} \, ;
\end{cases}
\end{gather*}
$S_{3}=\{(x^{3},u^{3})\}$ with
\begin{gather*}
u^{3}=
\begin{cases}
-1 & \text{ if } 0\leq t\leq 1 \, , \\
+1 & \text{ if } 1\leq t \leq T = 2 \, ;
\end{cases}\\
x^{3}_{1}=
\begin{cases}
-\frac{t^{2}}{2}+1 & \text{ if } 0\leq t\leq 1 \, , \\
\frac{t^{2}}{2}-2t+2 & \text{ if } 1\leq t\leq 2 \, ,
\end{cases}
\quad
x^{2}_{2}=
\begin{cases}
-t & \text{ if } 0\leq t\leq 1 \, ,\\
t-2 & \text{ if } 1\leq t\leq 2 \, .
\end{cases}
\end{gather*}
Direct calculations show that
\begin{equation*}
    I^{2}[x^{1},u^{1}]=[\sqrt{2},\sqrt{2}]=A \, ,
\end{equation*}
\begin{equation*}
    I^{2}[x^{2},u^{2}]=[4+\sqrt{6},-\sqrt{6}]=B \, ,
\end{equation*}
\begin{equation*}
    I^{2}[x^{3},u^{3}]=[2,0]=C \, .
\end{equation*}
Let $\zeta$ denote the set $S^{2}_{E}$. It is the continuous, convex
curve $\widehat{AB}$ (details can be found in \cite[\S
4.3]{Salukvadze}). As $C\in \zeta$ we have $S^{2}_{E} \cap S_{3}\neq
\emptyset$. Moreover, let us notice that $I_{3}$ has a form
$I_{3}[x,u]=f(I^{2}[x,u])$, where $f : \mathbb{R}^2 \rightarrow
\mathbb{R}$ is nondecreasing function. Therefore, the functional
$I_{3}$ is nonessential by Theorem \ref{t2}.

\begin{Remark}
If we change the functional $I_{3}$ into
\begin{equation}
\label{11}
  I_{3}[x,u]=\gamma_{1}I_{1}[x,u]+\gamma_{2}I_{2}[x,u],
\end{equation}
where $\gamma_{i}\in \mathbb{R}$ and $\gamma_{i}\geq 0, i=1,2$ or
\begin{equation}
\label{12}
 I_{3}[x,u]=[(I_{1}[x,u]-\sqrt{2})^{p}+(I_{2}[x,u]+\sqrt{6})^{p}]^{\frac{1}{p}},
\end{equation}
where $p\in [1,\infty]$, then  again $I_{3}$ will be nonessential by
Theorem \ref{t2}. It is worth noting that functionals (\ref{11}) and
(\ref{12}) can be used in order to find efficient solutions of the
problem $\min I^{2}[x,u]$ subject to \eqref{ex:cs}-\eqref{ex:ic}. We
mentioned this in section 2, details can be found in
\cite{Salukvadze} and \cite{Leitmann}.
\end{Remark}


\section{Conclusions}

The problem of optimizing a vector-valued criteria often arises in
connection with the solution of problems in the areas of planning,
organization of production, operational research and dynamical
control systems. Currently, the problem of optimizing a
vector-valued criteria is a central part of control theory and
great attention is being given to it in the design and
construction of modern automatic control systems, such as in
concrete applications of seismology, energetic chemistry and
metallurgy. In this work we use the notion of Pareto-optimality
in control theory to define and investigate nonessential
objective functionals of optimal control. For multicriteria
optimal control systems this notion seems to be new and not used
before. We claim the concept of nonessential objective functional
to be an important issue in optimal control and we trust it will
have an important role in the study of vector-valued optimization
problems of control theory. For future work, it would be
interesting to study the consequences of dropping nonessential
objectives in multi-criteria optimal control systems.


\section*{Acknowledgments}

Agnieszka B. Malinowska was supported by KBN under
Bia\l ystok Technical University grant No W/IMF/2/06;
Delfim F. M. Torres by the R\&D unit
``Centre for Research in Optimization and Control'' (CEOC).



\end{document}